\magnification=\magstep1
\baselineskip=14pt

\hoffset=-0.2cm

\centerline {\bf Two results on equations of nilpotent orbits.}
\bigskip
\centerline {by}
\bigskip
\centerline{\bf Jerzy Weyman%
\footnote*{\rm{supported by NSF, grant DMS 9700884.}}}
\bigskip

\beginsection 0. Introduction.

In this note I prove two results on the minimal generators of the defining
ideals of closure of nilpotent
conjugacy classes in the Lie algebra of $n\times n$ matrices. These results
were motivated by questions of Pappas
and Rapoport [PR] and they provide the answers to some of their conjectures.

In order to formulate the results let us recall the
following notation.

Let $K$ be a field and let $E$ be a vector space of dimension $n$ over
$K$. We denote by $X= Hom_K (E,E)$ the set of endomorphisms of $E$. Let
$\lbrace e_1 ,\ldots ,e_n\rbrace$ be a basis in $E$. We write the elements
of $X$ in the form $$\phi (e_j )=\sum_{i=1}^n \phi_{i,j}e_i ,$$ thus
identifying $X$ with a set of $n\times n$ matrices. We also identify the
space $Hom_K (E, E)$ with $E^*\otimes E$. In this tensor product we
identify $\phi$ with $$\phi =\sum_{i,j=1}^n \phi_{i,j} e_j^*\otimes e_i
.$$ The general linear group $GL(E)$ acts by conjugation on $X$. For a
partition $\mu =(\mu_1 ,\ldots ,\mu_s )$ of $n$ let $X_\mu$ be the closure
of the set $O(\mu )$ of nilpotent matrices with Jordan blocks of sizes
$\mu_1 ,\ldots ,\mu_s $. We denote by $A=K[X] = Sym (E\otimes E^* )$ the
coordinate ring of $X$, and $J(\mu )$ denotes the defining ideal of $X_\mu
$. The coordinate functions $\Phi_{i,j}$ are the generators of $A$. They
can be treated as independent indeterminates. We denote by $$\Phi =\left
(\matrix{\Phi_{i,j}}\right )_{1\le i,j\le n} $$ the generic $n\times n$
matrix.

Let us describe some important polynomials in $A$. Let $I,J$ be two
subsets of $[1,n]$ of cardinality $r$. Then $M(I,J)$ denotes the
determinant of the $r\times r$ submatrix of $\Phi$ with rows from $I$ and
columns from $J$.

We recall that the characteristic polynomial of $\phi$ can be written
$$P(\phi ,x)=\sum_{i=0}^{n} T_i x^{n-i}\ \ ,$$ where $T_0 =1$ and $T_i$ is
the sum of principal $i\times i$ minors of the generic matrix $\Phi$,
$$T_i = \sum_{I\subset[1,n], card(I)=i} M(I,I).$$
 The polynomials
$T_1 ,\ldots ,T_n$ are the generators of the ring of $GL(E)$-invariants in
$A$.

We are interested in  special nilpotent orbits. Let us fix a number $e<n$
and let $n=re+f$ be the division of
$n$ by $e$ with the remainder $f$. To such $n,e$ we associate the partition
$\mu (n,e) = (e^r ,f )$.
It is clear that the set $X_{\mu (n,e)}$ can be defined set theoretically
by the conditions
$\phi^e =0$. The subject of this note is the description of minimal
generators of the ideal $J(\mu (n,e))$.

\vskip .2cm

\proclaim Theorem 1. Let $K$ be a field of characteristic zero. Then the
ideal $J(\mu
(n,e))$ is generated  by the invariants $T_1 ,\ldots ,T_{e-1}$ and by the
entries of the matrix $\Phi^e$.
These polynomials form a minimal set of generators of $J(\mu (n,e))$.

\vskip .2cm

\proclaim Theorem 2. Let $K$ be an arbitrary  field and let $e=2$. The
ideal $J(\mu (n,2))$ is generated by the
invariants $T_1 ,T_2 ,\ldots ,T_n$ and by the entries of the matrix $\Phi^2$.

The proof of the first theorem is based on the geometric method used in
[W]. The main result in [W] described a small but nonminimal set of
generators of $J(\mu )$ for arbitrary $\mu$. Applying the method to the
smaller class of partitions $\mu (n,e)$ yields Theorem 1.

The proof of Theorem 2 relies on the paper [S] of Strickland. Using that
result one can easily show that the
ideal $J(n,2)$ is generated by the polynomials described in Theorem 2 and
by the minors of $\Phi$ of rank
$[{n\over 2}]+1$. It remains to show that the ideal generated by these
minors is contained in the ideal generated
by the polynomials described in Theorem 2. This is done by explicit
calculation, using some information about
 two-rowed Schur functors in characteristic $p>0$.

Throughout the paper we denote by $S_\mu E$ the Schur module corresponding
to the highest weight $\mu$ for the group $GL(n)$. This means that $\mu =
(\mu_1 ,\ldots ,\mu_n )$ where $\mu_i\in {\bf Z}$, $\mu_1 \ge\ldots\ge
\mu_n$.

\beginsection 1. The proof of Theorem 1.

Let $E, K$ be as above. For a given number $p$, $1\le p\le n$ we consider
the subrepresentation of $A$ which is the span of $p\times p$ minors of
the generic matrix $\Phi$. This vector space can be identified with
$\bigwedge^p E\otimes\bigwedge^p E^*$. For arbitrary $i$, $0\le i\le p$ we
define the $GL(E)$-equivariant map $\theta(i,p)$ to be the composition
$$\bigwedge^{i}E\otimes\bigwedge^{i}E^*\rightarrow
\bigwedge^{i}E\otimes\bigwedge^{i}E^*\otimes K\buildrel{1\otimes
tr^{(p-i)}}\over\rightarrow\bigwedge^{i}E\otimes\bigwedge^{i}E^*\otimes
\bigwedge^{p-i}E\otimes\bigwedge^{p-i}E^*\rightarrow
\bigwedge^{p}E\otimes\bigwedge^{p}E^* $$ where the second map is tensoring
with the $(p-i)$'st divided power of the trace element $tr=\sum_{i=1}^n
e_i\otimes e_i^*\in E\otimes E^*$ and the third map is a product of two
exterior multiplications.

We define $V_{i,p}$ to be the image of this map. When $K$ is a field of
characteristic zero it is well known that the space $\bigwedge^p
E\otimes\bigwedge^p E^*$ has the following $GL(E)$-equivariant
decomposition into irreducible representations of $GL(E)$, $$\bigwedge^p
E\otimes\bigwedge^p E^* =\oplus_{j=0}^{min(p,n-p)} S_{1^j ,0^{n-2j},
(-1)^j }E\ \ ,$$ where $S_{1^j ,0^{n-2j}, (-1)^j }E$ denotes the
irreducible representation (Schur functor) with the highest weight $(1^j
,0^{n-2j}, (-1)^j )$. We denote the summand $S_{1^j ,0^{n-2j}, (-1)^j }E$
occurring in $\bigwedge^p E\otimes\bigwedge^p E^*$ by $U_{j,p}$.

Notice that $U_{0,p}=V_{0,p}$ can be identified with the one dimensional
subspace generated by the invariant $T_p$.

\proclaim Lemma 1. Let $K$ be a field of characteristic zero. We have
$V_{i,p}=\sum_{j=0}^{min(i,n-p)}U_{j,p}$.
\par
\noindent{\sl Proof.} This is an easy calculation after one observes that
the highest weight vector in $\bigwedge^p E\otimes\bigwedge^p E^*$ of
weight $(1^i ,0^{2n-i}, (-1)^i )$ is $$w(i,p)=\sum_{J\subset [1,n],
|J|=p-i} (e_1\wedge\ldots\wedge e_i\wedge e_J)\otimes (e_J^*\wedge
e_{n+1-i}^*\wedge\ldots\wedge e_n^* ),$$ where for $J = \lbrace j_1
,\ldots ,j_t\rbrace$ we write $e_J =e_{j_1}\wedge\ldots\wedge e_{j_t}$,
$e_J^* =e_{j_1}^*\wedge\ldots\wedge e_{j_t}^*$.$\bullet$

\proclaim Remark.
\item{a)} Notice that $U_{0,p}=V_{0,p}$ can be identified with the one
dimensional subspace generated by the invariant $T_p$. The representation
$V_{1,p}$ is the subspace generated by the entries of the $p$-th power
$\Phi^p$ of the generic matrix.
\item{b)} We have a commutative diagram
$$\matrix{\bigwedge^{i-1}E\otimes\bigwedge^{i-1}E^*&\buildrel{1\wedge
T_1}\over\longrightarrow&\bigwedge^i E\otimes\bigwedge^i E^*\cr\cr
\llap{$\scriptstyle (p-i+1)\theta(i-1,p)$}
\searrow&&\swarrow\rlap{$\scriptstyle  \theta(i,p)$}\cr\cr &\bigwedge^p
E\otimes\bigwedge^p E^*&}$$ which proves that $V_{i-1,p}\subset V_{i,p}$.

Before we start with the proof we need a couple of elementary lemmas.

\proclaim Lemma 2. Let $i\ge 1$. Then the vector space $V_{i,p+1}$ is
contained in the ideal generated by
$V_{i,p}$.
\par
\noindent{\sl Proof.} Let us take a typical element of $V_{i,p+1}$. It can
be written as $$\sum_{j_1 ,\ldots ,j_{p+1-i}}M(u_1 ,\ldots u_i ,j_1
,\ldots ,j_{p+1-i}; v_1 ,\ldots ,v_i ,j_1 ,\ldots ,j_{p+1-i}).$$ Here we
are summing over all choices of $j_1 ,\ldots ,j_{p+1-i}$ regardless of
whether the summands are nonzero. Let us take the Laplace expansion of
every term in the above sum, with respect to the row $u_1$. This allows us
to express our generator in the following form $$\sum_{j_1 ,\ldots
,j_{p+1-i}}(\sum_{k=1}^i (-1)^{k+1}\Phi_{u_1 ,v_k}M(u_2 ,\ldots u_i ,j_1
,\ldots ,j_{p+1-i}; v_1 ,\ldots ,{\hat v_k}\ldots ,v_i ,j_1 ,\ldots
,j_{p+1-i})+$$ $$+\sum_{l=1}^{p+1-i}(-1)^{i+l+1}\Phi_{u_1 ,j_l} M(u_2
,\ldots u_i ,j_1 ,\ldots ,j_{p+1-i}; v_1 ,\ldots ,v_i ,j_1 ,\ldots ,{\hat
j}_l ,\ldots ,j_{p+1-i})).$$

Writing the first part of the sum as a linear combination  of $\Phi_{u_1
,v_k}$ we see that all coefficients
are in $V_{i-1,p}\subset V_{i,p}$.  Similarly, writing the second part of
the sum as a linear combination of
$\Phi_{u_1 ,j_l}$ we see that each coefficient is in $V_{i,p}$.$\bullet$
\medskip
We also recall the geometric method for calculating syzygies. We consider
the Springer type desingularization $Y_\mu$ of $X_\mu$. Let $\mu^\prime
=(\mu^\prime_1 ,\ldots ,\mu^\prime_t )$ be the partition conjugate to
$\mu$. Let $G/P_{\mu^\prime}$ be the flag variety of partial flags $(R_1
\subset\ldots\subset R_{t-1}\subset R_t =E)$ with $dim\ R_i =\mu^\prime_1
+\ldots +\mu^\prime_i$, and $$Y_\mu =\lbrace (\phi ,R_1 ,\ldots ,R_t)\in
X\times G/P_{\mu^\prime}\ |\ \phi (R_i )\subset R_{i-1}\ {\rm for}\
i=1,\ldots ,t\ \rbrace .$$

We denote by ${\cal R}_i$ the tautological vector bundle on
$G/P_{\mu^\prime}$ of dimension $\mu^\prime_1 +\ldots
+\mu^\prime_i$, corresponding to the subspace $R_i$. We denote by ${\cal
Q}_i$ the tautological factorbundle
$(E\times G/P_{\mu^\prime})/{\cal R}_i$.

The variety $Y_\mu$ is the total space of the cotangent bundle on
$G/P_{\mu^\prime}$. Let ${\cal T}_\mu$ denote the tangent bundle on
$G/P_\mu$. Then we have an exact sequence of vector bundles on
$G/P_{\mu^\prime}$ $$0\rightarrow {\cal S}_\mu\rightarrow E\otimes
E^*\times G/P_{\mu^\prime}\rightarrow {\cal T}_\mu\rightarrow 0.$$

The bundle ${\cal S}_\mu$ is ``the staircase bundle". In terms of
tautological bundles ${\cal S}_\mu$ can be described (as a subbundle of
$(E\otimes E^* )\times G/P_{\mu^\prime}$) as $${\cal S}_\mu = {\cal
R}_1\otimes {\cal Q}_0^* +{\cal R}_2\otimes {\cal Q}_1^* +\ldots + {\cal
R}_t \otimes {\cal Q}_{t-1}^*$$ (the sum is not direct), with the
convention ${\cal R}_t = E\times G/P_{\mu^\prime}$, ${\cal Q}_0 =
E^*\times G/P_{\mu^\prime}$.

Recall that the basic theorem of the geometric method (comp. [W], [F-W])
implies the following

\proclaim Theorem 3. The varieties $X_\mu$ are normal, with rational
singularities. Moreover, a minimal free
resolution of the coordinate ring $K[X_\mu ]$ as an $A$-module has the terms
$$\ldots \rightarrow F^\mu_i\rightarrow F^\mu_{i-1}\rightarrow\ldots
\rightarrow F^\mu_1\rightarrow
F^\mu_0\rightarrow 0$$ where
$$F^\mu_i =\oplus_{j\ge 0} H^j (G/P_\mu ,\bigwedge^{i+j}{\cal S}_\mu
)\otimes_K A(-i-j).$$
where $A(-d)$ denotes the one dimensional homogeneous free $A$-module with
a generator in degree $d$.
\par
Let ${\hat\mu}$ be the partition ${\hat\mu}=(max(\mu_1 -1,0),max(\mu_2
-1,0),\ldots ,max(\mu_t -1,0))$. Then ${\hat\mu}^\prime = (\mu_2^\prime
,\ldots ,\mu_t^\prime )$. We see that we have an exact sequence
$$0\rightarrow E^*\otimes {\cal R}_1\rightarrow {\cal S}_\mu \rightarrow
{\cal S}_{\hat\mu}({\cal Q}_1 ,{\cal Q}_1^* )\rightarrow 0\eqno{(*)}$$
where ${\cal S}_{\hat\mu}({\cal Q}_1 ,{\cal Q}_1^* )$ is the bundle ${\cal
S}_{\hat\mu}$ constructed in a relative situation with ${\cal Q}_1$
replacing $E$.

Using the exact sequence $(*)$ we can estimate the terms of the complex
$F^\mu_\bullet$ inductively. We consider
the minimal resolution $F^{\hat\mu}({\cal Q}_1 ,{\cal Q}_1^* )_\bullet$ of
the nilpotent orbit closure for the
partition
$\hat\mu$ in the relative situation (i.e. ${\cal Q}_1$ replaces $E$ and we
treat the resolution as the sequence
of $Sym ({\cal Q}_1\otimes {\cal Q}_1^* )$-modules). For each term
$S_\lambda {\cal Q}_1$ occurring in that resolution we can describe the
push down $K_\lambda (E, E^*)$ of the
complex $S_{\lambda}{\cal Q}_1\otimes\bigwedge^\bullet (E^*\otimes {\cal
R}_1 )$. The terms of the complexes
$K_\lambda (E,E^*)$ for all representations $S_\lambda {\cal Q}_1 $
occurring in $F^{\hat\mu}({\cal Q}_1 ,{\cal
Q}_1^* )_\bullet$ give an upper bound on the terms from $F^\mu_\bullet$.

The following facts are proved in [W], section 3.

\proclaim Theorem 4. Let $\mu$, $\hat\mu$ be as above.
 Let $S_\lambda {\cal Q}_1$ be the term in $F^{\hat\mu}({\cal Q}_1 ,{\cal
Q}_1^* )_i$. Then the terms
coming from $K_\lambda (E, E^* )$ appear in homological degrees $\ge i$,
i.e. they can give the contribution only
to the terms $F^\mu_j$ with $j\ge i$.
\par
In particular, as stated in Corollary (3.15) of [W], the term $F^\mu_0$ is
trivial, consisting of a copy of $A$ in homogeneous degree $0$. The terms
in $F^\mu_1$ are, by Theorem (3.12) of [W],  of two types:
\item {1.} The term  $F^{\hat\mu}({\cal Q}_1 ,{\cal Q}_1^* )_0$ is trivial,
so the corresponding complex $K_{(0)}$
is just a pushdown of
$\bigwedge^\bullet (E^*\otimes {\cal R}_1 )$. Its first term gives a
possible contribution to $F^\mu_1$
which is the vanishing of
$n-\mu_1 +1$ minors of the matrix $\Phi$,
\item {2.} Each term $S_\lambda {\cal Q}_1$ from $F^{\hat\mu}({\cal Q}_1
,{\cal Q}_1^* )_1$ gives a possible contribution to $F^\mu_1$ which is the
zero'th term of the corresponding complex $K_\lambda$.
 \medskip Now let us
use the above inductive procedure for the family of partitions $\mu (n,e)
= (e^r ,f )$. The partition ${\hat\mu} (n,e) = ((e-1)^r ,max(f-1,0))$ is
just $\mu (n-r-1,e-1)$ for $f>0$ or $\mu (n-r,e-1)$ for $f=0$. These
partitions are in the same family, so we can assume by induction that
Theorem 1 is true for ${\hat\mu}$.

The inductive assumption means that the first term of $F^{\hat\mu
(n,e)}_1$ consists of the trivial representations in homogeneous degrees
$1,\ldots ,e-1$ and the representation $S_{1,0,\ldots ,0,-1}{\cal Q}_1$ in
homogeneous degree $e-1$.  This means that the terms in $F^\mu_1$ have
trivial representations in homogeneous degrees $1,\ldots ,e-1$, the
representation $E^*\otimes E$ in degree $e$ (comp. Lemma (3.11) in [W])
and the terms coming from $n-\mu(n,e)_1 +1=n-r$ size minors of the matrix
$\Phi$. The generators in degrees $\le e$ have to consist of the
invariants $T_1 ,\ldots ,T_{e-1}$ and the vanishing of entries of
$\Phi^e$, because these terms have to occur in $F^\mu_1$ and they match
the terms we described. Thus the induction implies that the minimal
generators of $J(\mu (n,e))$ are those listed in Theorem 1 plus possibly
some linear combinations of $(n-r)\times (n-r)$ minors of $\Phi$.

However we also recall what the Theorem (4.6) from  [W] says about the
ideal $J(n,r)$.

\proclaim Theorem 5. The ideal $J(\mu (n,e))$ is generated (nonminimally)
by the invariants $U_{0,p}$ ($1\le p\le
n$) and by the representations $U_{i,ie-i+1}$ (for $1\le i\le r$).
\par
The only possible generator in degree $n-r$ on the above list is the
invariant $U_{0,n-r}$. However this cannot be a minimal generator because
$U_{0,n-r}$ is contained in $V_{1,n-r}$ which is in the ideal generated by
$V_{1,e}$ according to Lemma 2. This concludes the proof of Theorem 1.

\beginsection 2. The proof of Theorem 2.

In this section we work over a field $K$ of arbitrary characteristic.
We are interested in the generators of the ideal $J(\mu (n,2))$. Let us
recall that Strickland in [S] exhibited a
standard basis in the coordinate ring $A/J(\mu (n,2))$ of the orbit closure
$X_{\mu (n,2)}$.
Denote by $M(a_1 ,\ldots ,a_r\ ;b_1 ,\ldots ,b_r )$ the $r\times r$ minor
of $\Phi$ corresponding to rows $a_1
,\ldots ,a_r$ and columns $b_1 ,\ldots ,b_r $. One notices easily that the
only relations used in [S] to prove that
the standard monomials span the ring
$A/J(\mu (n,2))$ are the relations

$$\sum_{1\le i_1 < \ldots <i_r \le n}M(a_1 ,\ldots ,a_{p-r},i_1 ,\ldots
,i_r\ ; b_1 ,\ldots ,b_{p-r},i_1 ,\ldots ,i_r )\eqno{Rel(r,p)}$$ for $1\le
p\le n, 1\le r\le p$ and $([{n\over 2}]+1)\times ([{n\over 2}]+1)$ minors
of $\Phi$.  Therefore Theorem 2 is a consequence of the following lemmas.

\proclaim Lemma 3. The polynomials of type $Rel(r,p)$ for $1\le p\le n,1\le
r\le p$ are in the ideal generated by
the invariants $T_1 ,T_2 ,\ldots ,T_n$ and by the entries of $\Phi^2$.
\par
\noindent{\sl Proof.} Let us denote by $J$ the ideal generated by the
invariants $T_1 ,\ldots ,T_n$ and by the entries of $\Phi^2$. We use
induction on $p$. For $p=1,2$ we see that the equations of type $Rel(r,p)$
describe exactly $T_1$ and the entries of $\Phi^2$. Let us assume that the
statement is true for $p$. Let us first assume that $r>1$. If $p-r=0$ we
are dealing with the invariant $T_p$, so it is in $J$. If $p-r>0$ we see
that after taking the Laplace expansion of $Rel(r,p+1)$ with respect to
the row $a_1$ we get a linear combination of terms of type $Rel(r-1,p)$ so
we are done by induction. For $r=1$ we note that

$$\sum_{1\le i\le n} M(a_1 ,\ldots ,a_p ,i;b_1 ,\ldots ,b_p ,i )=$$
$$=\sum_{1\le i\le n}\sum_{1\le u\le p} (-1)^{u+1}M(a_1 ;b_u )M(a_2 ,\ldots
,a_p ,i;b_1 ,\ldots ,{\hat b}_u
,\ldots ,b_p ,i)+$$
$$+(-1)^p \sum_{1\le i\le n} M(a_1 ;i)M(a_2 ,\ldots ,a_p ,i;b_1 ,\ldots
,b_p ).$$
The first summand on the right hand side is in $J$ by induction. It remains
to deal with the second summand which
we will denote by $F$. We take the Laplace expansion with the respect to
the row $i$ in each term in that summand.
We obtain
$$F=\sum_{1\le i\le n}\sum_{1\le v\le p} (-1)^{v+1} M(a_1 ;i)M(i;b_v )M(a_2
,\ldots ,a_p ;b_1 ,\ldots ,{\hat b}_v
,\ldots ,b_p ).$$
However
$$\sum_{1\le i\le n}M(a_1 ;i)M(i;b_v ) =  \sum_{1\le i\le n}M(i
;i)M(a_1;b_v )+\sum_{1\le i\le n} M(a_1 ,i;b_v
,i)\ \ .$$ Since the second summand above is in $J$, we see that in $A/J$
we have $$F= \sum_{1\le i\le n} M(i,i)\sum_{1\le v\le p}(-1)^{v+1} M(a_1
,b_v )M(a_2 ,\ldots ,a_p ;b_1 ,\ldots ,{\hat b}_v ,\ldots ,b_p )=$$ $$=
T_1 M(a_1 ,\ldots ,a_p ;b_1 ,\ldots ,b_p )$$ so $F\in J$ and we are done.

\proclaim Lemma 4. The $([{n\over 2}]+1)\times ([{n\over 2}]+1)$ minors of
$\Phi$ are contained in the ideal
spanned by the invariants $T_1 ,T_2 ,\ldots ,T_n$ and by the entries of the
matrix $\Phi^2$.
\par
\noindent{\sl Proof.} We employ the following notation. We write $m=
[{n\over 2}]+1$. We identify the set of $m\times m$ minors with
$\bigwedge^m E\otimes\bigwedge^m E^* =\bigwedge^m
E\otimes\bigwedge^{n-m}E$. This is only an $SL(E)$-isomorphism but it is
sufficient for our purposes. The relation $Rel(r,m)$ translates under this
identification to the image of the following map.
$$\psi(r,m):\bigwedge^{m-r}E\otimes\bigwedge^{n-m+r}E\buildrel{1\otimes\Delta}\over\rightarrow
\bigwedge^{m-r}E\otimes\bigwedge^r
E\otimes\bigwedge^{n-m}E\buildrel{m\otimes 1}\over\rightarrow \bigwedge^m
E\otimes\bigwedge^{n-m}E$$

Therefore we need to prove

\proclaim Lemma 5. The vector space $\bigwedge^m E\otimes\bigwedge^{n-m}E$
is spanned by the images of the maps $\psi (r,m)$ for $1\le r\le m$.

\noindent{\sl Proof.} We use induction on $n$. We prove the spanning
weight by weight. Notice that the only weights in $\bigwedge^m
E\otimes\bigwedge^{n-m}E$ are (up to permutation of the basis) the weights
$(2^j ,1^{n-2j},0^j )$. Notice that the problem of spanning in the weight
$(2^j ,1^{n-2j},0^j )$ for a given $n$ is equivalent to the problem of
spanning in the weight $(1^{n-2j})$ for the dimension $n-2j$. So by
induction on $n$ we can assume that all vectors of weights $\ne (1^n )$
are in the subspace generated by the images of $\psi(r,m)$. But the
subrepresentation of $\bigwedge^m E\otimes\bigwedge^{n-m}E$  generated by
all vectors of weights $\ne (1^n )$ is isomorphic to the kernel of the
exterior multiplication $$\bigwedge^m
E\otimes\bigwedge^{n-m}E\rightarrow\bigwedge^n E.$$ Therefore it is enough
to show that the images of $\psi(r,m)$ cannot be all contained in this
kernel. However one sees easily that the map $\psi (r,m)$ composed with
the exterior multiplication is the exterior multiplication
$$\bigwedge^{m-r}E\otimes\bigwedge^{n-m+r}E\rightarrow\bigwedge^n E$$
multiplied by the binomial coefficient ${n-m+r}\choose r$. So Theorem 2
follows from the last elementary lemma.

\proclaim Lemma 6. Let $n$ be an integer, $m=[{n\over 2}]+1$. Then
$$GCD_{1\le r\le m}{{n-m+r}\choose r}=1.$$

\noindent{\sl Proof.} Let $n=2t$ be even. Then $m=t+1$. We are interested
in $$GCD_{1\le r\le t+1}{{t-1+r}\choose r}.$$ Using the relations
${k\choose r}={{k-1}\choose r}+{{k-1}\choose {r-1}}$ we get $$GCD_{2\le
r\le t+1}{{t-1+r}\choose r}= GCD_{2\le r\le t+1} {{t+1}\choose r}$$ which
obviously equals to $1$.

If $n=2t+1$ is odd, then $m=t+1$ and we are interested in $$GCD_{1\le r\le
t+1}{{t+r}\choose r}.$$ Again using the same method we find that
$$GCD_{1\le r\le t+1}{{t+r}\choose r}=GCD_{1\le r\le t+1}{{t+1}\choose r}
=1.$$

\beginsection References.

\item{[FW]}  T.~Fukui, J.~Weyman, Cohen-Macaulay properties of
Thom-Boardman strata. I. Morin's ideal.
Proc. London Math. Soc. (3) 80  (2000), no. 2, 257-303,

\item{[J]} G.~James, Representations of the Symmetric Group, Springer LN, 628,

\item{[PR]} G.~Pappas, M.~Rapoport, Local models in the ramified case. I.
The EL-case, preprint, 2000,

\item{[S]} E.~Strickland, On the variety of projectors, J. Algebra,
106(1987), 135-147,

\item{[W]}  J.~Weyman, The equations of conjugacy classes of nilpotent
matrices ,
       Invent. Math. 98, 229-245 (1989),

\vskip 3cm
\leftline{Department of Mathematics}
\leftline{Northeastern University}
\leftline{Boston, MA, 02115}
\leftline{USA}
\leftline{e-mail:weyman@research.neu.edu}

\bye